\documentclass[11pt, twoside]{article}
\usepackage{latexsym}
\usepackage{amsmath}
\usepackage{amssymb}
\usepackage[all]{xy}
\usepackage{amsfonts}
\usepackage{verbatim}
\usepackage{amsthm}
\usepackage{mathrsfs}
\usepackage{epsfig}
\usepackage{xy}
\usepackage{array}
\usepackage{stmaryrd}
\usepackage{graphicx,color}
\usepackage{xcolor}
\usepackage{tikz}
\usetikzlibrary{arrows,calc}
\usepackage{etex}
\usepackage{mathdots}
\usepackage{float}
\usepackage{graphics}
\usepackage{pdflscape}
\usepackage{xcolor}
\usepackage[colorlinks=true,linkcolor=blue,citecolor=blue]{hyperref}
\usepackage{anysize,hyperref}
\input xypic
\xyoption{all}

\usepackage[perpage,symbol]{footmisc}
\usepackage{setspace}
\def\A{\mathcal{A}}

\def\B{\mathcal{B}}
\def\M{\mathcal{M}}

\def\C{\mathscr{C}}
\def\E{\mathbb{E}}
\def\s{\mathfrak{s}}

\def\op{^\mathrm{op}}

\def\dr{\ar@{->}[r]}

\def\X{\mathcal{X}}
\def\Y{\mathcal{Y}}
\def\Z{\mathcal{Z}}
\def\T{\mathcal{T}}

\newcommand{\ppr}{^{\prime}}          %%% '
\newcommand{\tri}[7]{\xymatrix@C=1.5em{#1\ar[r]^{#5}&#2\ar[r]^{#6}&#3\ar[r]^{#7}&#4}}
\newcommand{\ltri}[8]{\xymatrix@C=#8cm{#1\ar[r]^{#5}&#2\ar[r]^{#6}&#3\ar[r]^{#7}&#4}}
\newcommand{\rtri}[6]{\xymatrix@C=1.5em{#1\ar[r]^{#4}&#2\ar[r]^{#5}&#3\ar[r]^{#6}&\Sigma #1}}

\begin{document}
\baselineskip=15pt
\title{\Large{\bf Relative tilting theory in extriangulated categories}}
\medskip
\author{Chenbei Xie}

\date{}

\maketitle
\def\blue{\color{blue}}
\def\red{\color{red}}

\newtheorem{theorem}{Theorem}[section]
\newtheorem{lemma}[theorem]{Lemma}
\newtheorem{corollary}[theorem]{Corollary}
\newtheorem{proposition}[theorem]{Proposition}
\newtheorem{conjecture}{Conjecture}
\theoremstyle{definition}
\newtheorem{definition}[theorem]{Definition}
\newtheorem{question}[theorem]{Question}
\newtheorem{remark}[theorem]{Remark}
\newtheorem{remark*}[]{Remark}
\newtheorem{example}[theorem]{Example}
\newtheorem{example*}[]{Example}
\newtheorem{condition}[theorem]{Condition}
\newtheorem{condition*}[]{Condition}
\newtheorem{construction}[theorem]{Construction}
\newtheorem{construction*}[]{Construction}

\newtheorem{assumption}[theorem]{Assumption}
\newtheorem{assumption*}[]{Assumption}

\baselineskip=17pt
\parindent=0.5cm
\vspace{-6mm}

\begin{abstract}
\baselineskip=16pt
In this article, we define relative resolutions and coresolutions in  extriangulated categories.
By studying this relative resolutions and coresolutions, we get a generalization of the Auslander-Buchweitz approximation theory. Finally, we develop some theories of relative tilting objects in extriangulated categories. \\[1mm]
\end{abstract}

\pagestyle{myheadings}
\markboth{\rightline {\scriptsize   Chenbei Xie }}
         {\leftline{\scriptsize Relative tilting theory in extriangulated categories}}

\section{Introduction}
Exact categories (abelian categories are also exact categories) and triangulated categories are two fundamental structures in algebra and geometry. By extracting those properties of ${\rm Ext}^1(-,-)$ on exact categories and triangulated categories,  Nakaoka and Palu \cite{NP} introduced the notion of extriangulated categories, which is a simultaneous generalization of exact categories and triangulated categories.

The approximation theory has its origin with the concept of injective envelopes and
it has had a wide development in the context of module categories.
Auslander and Buchweitz \cite{AB} studied the ideas of injective envelopes and
projective covers in terms of maximal Cohen-Macaulay approximations for certain
modules. Based on Auslander-Buchweitz's idea, Hashimoto \cite{H} defined the so called ``Auslander-Buchweitz context" for abelian categories, and provided a new framework to homological approximation theory.
In \cite{MSSS}, they introduced and  developed an analogous of the Auslander-Buchweitz
approximation theory in the context of triangulated categories by using a version
of relative homology.

Recently,  Monroy and Hern\'{a}ndez \cite{MH} introduced a special kind of relative resolutions and
coresolutions associated to a pair of classes of objects in an abelian category.
By studying this resolutions and coresolutions they has obtained a generalization of the Auslander-Buchweitz approximation theory. They also showed that a generalization of Salce's Lemma and Garcia-Rozas's Lemma, among other classical results from Auslander-Buchweitz and Auslander-Reiten theory.
In this article, based on the work of Monroy and Hern\'{a}ndez \cite{MH} we introduce and study a relative tilting theory in an extriangulated category. Specifically, assume that $\X$ and $\Y$ are two classes of objects in an extriangulated category. We introduce $(\X,\Y)$-resolutions (coresolutions) and resolutions (coresolutions) dimension. We use these concepts to generalize the Auslander-Buchweitz-Reiten approximations theory.

This article is organized as follows. In Section 2, we review some elementary definitions and facts on extriangulated categories. In Section 3, we give the definitions of relative homological dimensions,
the class ${\rm Fac}^{\X}_n(\T)$, $\X$-complete pairs and $\X$-hereditary pairs where $\X$ and $\T$ are two classes of an extriangulated category. Moreover, we also present the notion of the $n$-$\X$-tilting classes and study its good properties.

\section{Preliminaries}
Throughout this article,
we fix a commutative ring $R$ with $1$, and $(\C,\E,\s)$ denotes an
$R$-linear  extriangulated category
defined in \cite{NP} (see Section 2 of \cite{NP} for more details).
We will use the following terminology and known results in extriangulated category.

\begin{definition}
For any class of objects $\X\subseteq\C$ and any $i\geqslant 1$, we define the right $i$-th orthogonal class $\X^{\perp_i}:=\{C\in\C\mid\E^i(-,C)|_{\X}=0\}$ and the right orthogonal class $\X^{\perp}:=\bigcap_{i>0}\X^{\perp_i}$ of $\X$.
Dually, we have the left $i$-th orthogonal class $^{\perp_i}\X$ and the left orthogonal class $^{\perp}\X$ of $\X$.

With respect inclusion of classes and objects, $C\in\C$ means that C is an object of $\C$, $\M\subseteq\C$ means that $\M$ is a class of objects of $\C$, $(\A,\B)\subseteq \C^2$ means that $\A\subseteq \C$ and $\B\subseteq \C$.

For some $\M\subseteq\C$, we have some classes of objects in $\C$. The class smd($\M$) is the class whose objects are all the direct summands of objects in $\M$ , the class $\M^{\oplus}(\M^{\oplus <\infty})$ of all the (finite) coproducts of objects in $\M$ .
\end{definition}

\begin{remark}
Recall that  a right $\C$-module is a contravariant additive functor $\C^{op} \rightarrow {\rm Mod}R$ to the category of $R$-modules. Dually, a left $\C$-module is a covariant additive functor to $R$-modules.
All  right $\C$-modules and left $\C$-modules form abelian categories ${\rm Mod}\C$ and ${\rm Mod}\C^{op}$,
respectively. Moreover, these categories have enough projectives and they are given by direct summands of direct sums of representable functors $\C(-,A)={\rm Hom}_{\C}(-,A)\in {\rm Mod}\C$, $\C(A,-)={\rm Hom}_{\C}(A,-)\in \C{\rm Mod}$.
\end{remark}

\begin{theorem}
[\cite{GNP}, Chapter III, Theorem 3.5]Let $(\C, \mathbb{E}, \s)$ be an extriangulated category and $\xymatrix{A\ar[r]^{x}&B\ar[r]^{y}&C\ar@{-->}[r]^{\delta}&}$ be any $\mathbb{E}$-triangle in $\C$.

{\rm (1)} We have a long exact sequence
$$\C(-,A)\rightarrow\C(-,B)\rightarrow\C(-,C)\rightarrow\E(-,A)\rightarrow\cdots$$
$$\cdots\rightarrow\E^{n-1}(-,C)\rightarrow\E^{n}(-,A)\rightarrow\E^{n}(-,B)\rightarrow\E^{n}(-,C)\rightarrow\cdots$$
in Mod $\C$.

{\rm (2)} Dually, we have a long exact sequence
$$\C(C,-)\rightarrow\C(B,-)\rightarrow\C(A,-)\rightarrow\E(C,-)\rightarrow\cdots$$
$$\cdots\rightarrow\E^{n-1}(A,-)\rightarrow\E^{n}(C,-)\rightarrow\E^{n}(B,-)\rightarrow\E^{n}(A,-)\rightarrow\cdots$$
in $\C$Mod.

\end{theorem}

We give the following notions of the relative projective (injective) dimension.

\begin{definition}
Let $\C$ be an extriangulated category, $\B, \A\subseteq\C$ and $C\in\C$.
\begin{itemize}
\item[(1)] The $\A$-projective dimension of $C$ is ${\rm pd}_{\A}(C):= {\rm min}\{n\in\mathbb{N}:\mathbb{E}^{k}(C,-)|_{\A}=0,\forall k>n \}$;
the $\A$-projective dimension of $\B$ is ${\rm pd}_{\A}(\B):={\rm sup}\{{\rm pd}_{\A}(B):B\in\B\}$ .
The $\A$-injective dimension ${\rm id}_{\A}(C)$ of C and $\A$-injective dimension ${\rm id}_{\A}(\B)$ of $\B$ are defined dually.

\item[(2)] For a pair $(\X,\omega)\subseteq\C^2$, it is said that $\omega$ is a relative cogenerator in $\X$ if $\omega\subseteq\X$ and any $X\in\X$ admits an $\mathbb{E}$-triangle $\xymatrix{X\ar[r]&W\ar[r]&X'\ar@{-->}[r]^{\delta}&}$ in $\C$, with $W\in\omega$ and $X'\in\X$.
$\omega$ is $\X$-injective if ${\rm id}_{\X}(\omega)=0$.
Dually, we can define the notions of relative generator in $\X$ and $\X$-projective.
\end{itemize}

\end{definition}

\begin{definition}
For a class $\Z$ of objects in an extriangulated category $\C$, a morphism $f\colon Z\rightarrow M$ in $\C$ is called a $\Z$-precover if $Z\in\Z$ and ${\rm Hom}_{\C}(Z',f):{\rm Hom}_{\C}(Z',Z)\rightarrow {\rm Hom}_{\C}(Z',M)$ is an epimorphism $\forall Z'\in\Z$. And it is said that $\Z$ is precover if each $C\in\C$ admits a $\Z$-precover $Z\rightarrow C$.
The notions of $\Z$-preenvelop and preenvelop class are defined dually.
\end{definition}

The following are two important propositions that we often use in our main results.

\begin{proposition}\label{RemET4op}{\rm
$($\cite{NP}, Paraphrase of {\rm (ET4)$\op$}$)$}

Let $\xymatrix{D\ar[r]^{f'}&A\ar[r]^f&B\ar@{-->}[r]^{\delta}&}$ and $\xymatrix{F\ar[r]^{g'}&B\ar[r]^g&C\ar@{-->}[r]^{\delta'}&}$ be $\E$-triangles. Then there exist an $\E$-triangle $\xymatrix{E\ar[r]^{h'}&A\ar[r]^h&C\ar@{-->}[r]^{\delta''}&}$ and a commutative diagram
\[
\xy
(-21,7)*+{D}="0";
(-7,7)*+{E}="2";
(7,7)*+{F}="4";
(-21,-7)*+{D}="10";
(-7,-7)*+{A}="12";
(7,-7)*+{B}="14";
(-7,-21)*+{C}="22";
(7,-21)*+{C}="24";
{\ar^{d} "0";"2"};
{\ar^{e} "2";"4"};
{\ar@{=} "0";"10"};
{\ar_{h\ppr} "2";"12"};
{\ar^{g\ppr} "4";"14"};
{\ar_{f\ppr} "10";"12"};
{\ar_{f} "12";"14"};
{\ar_{h} "12";"22"};
{\ar^{g} "14";"24"};
{\ar@{=} "22";"24"};
{\ar@{}|\circlearrowright "0";"12"};
{\ar@{}|\circlearrowright "2";"14"};
{\ar@{}|\circlearrowright "12";"24"};
\endxy
\]
in $\C$, satisfying the following compatibilities.
\begin{itemize}
\item[{\rm (i)}] $\xymatrix{D\ar[r]^{d}&E\ar[r]^e&F\ar@{-->}[r]^{g'*\delta}&}$ is an $\E$-triangle,
\item[{\rm (ii)}] $\delta'=e_{*} \delta '' $,
\item[{\rm (iii)}] $d_{*}\delta=g^{*}\delta''$.
\end{itemize}
\end{proposition}

\begin{proposition}\label{prop2}{\rm (\cite{NP}, Proposition 3.15)}
Let $\C$ be an extriangulated category. Then the following holds.
Let $C$ be any object, and let
$\xymatrix{A_1\ar[r]^{x_1}&B_1\ar[r]^{y_1}&C\ar@{-->}[r]^{\delta_1}&}$, $\xymatrix{A_2\ar[r]^{x_2}&B_2\ar[r]^{y_2}&C\ar@{-->}[r]^{\delta_2}&}$
be any pair of $\E$-triangles. Then there is a commutative diagram in $\C$
\[
\xy
(-7,21)*+{A_2}="-12";
(7,21)*+{A_2}="-14";
(-21,7)*+{A_1}="0";
(-7,7)*+{M}="2";
(7,7)*+{B_2}="4";
(-21,-7)*+{A_1}="10";
(-7,-7)*+{B_1}="12";
(7,-7)*+{C}="14";
{\ar@{=} "-12";"-14"};
{\ar_{m_2} "-12";"2"};
{\ar^{x_2} "-14";"4"};
{\ar^{m_1} "0";"2"};
{\ar^{e_1} "2";"4"};
{\ar@{=} "0";"10"};
{\ar_{e_2} "2";"12"};
{\ar^{y_2} "4";"14"};
{\ar_{x_1} "10";"12"};
{\ar_{y_1} "12";"14"};
{\ar@{}|\circlearrowright "-12";"4"};
{\ar@{}|\circlearrowright "0";"12"};
{\ar@{}|\circlearrowright "2";"14"};
\endxy
\]
which satisfies
\begin{eqnarray*}
&\s(y_2^*\delta_1)=[\xymatrix{A_1\ar[r]^{m_1}&M\ar[r]^{e_1}&B_2}],&\\
&\s(y_1^*\delta_2)=[\xymatrix{A_2\ar[r]^{m_2}&M\ar[r]^{e_2}&B_1}],&\\
&m_{1*}\delta_1+m_{2*}\delta_2=0.&
\end{eqnarray*}

\end{proposition}

\section{Main results}
Let us introduce the definition of relative resolutions and coresolutions and we will use it for the main result to give a generalization of Auslander-Buchweitz and Auslander-Reiten theory.

\begin{definition}
Let $\C$ be an extriangulated category, $M\in\C, \X, \Y$ and  $ \Z \subseteq\C$.

If we have the following $\mathbb{E}$-triangles in $\C$:
$$\xymatrix{M\ar[r]&Y_0\ar[r]&N_0\ar@{-->}[r]&}$$
$$\xymatrix{N_0\ar[r]&Y_1\ar[r]&N_1\ar@{-->}[r]&}$$
$$\vdots$$
$$\xymatrix{N_{n-2}\ar[r]&Y_{n-1}\ar[r]&N_{n-1}\ar@{-->}[r]&}$$
$$\xymatrix{N_{n-1}\ar[r]&Y_{n}\ar[r]&N_{n}\ar@{-->}[r]&}$$
$$\vdots$$

(1) We can get a complex $M\rightarrow Y_0\rightarrow Y_1\rightarrow\cdots\rightarrow Y_{n-1}\rightarrow Y_n\rightarrow\cdots$ in $\C$, with $Y_k\in\Y\cup\{0\}$ and $N_k\in\X\cup\{0\} \forall k\geq 0$, is called a $\Y_{\X}$-coresolution of M.

(2) A complex $M\rightarrow Y_0\rightarrow Y_1\rightarrow\cdots\rightarrow Y_{n-1}\rightarrow Y_n$ in $\C$, with $Y_n\in\X\cap\Y, Y_k\in\Y$ and $N_k\in\X, \forall k\in[0,n-1]$, is called a finite $\Y_{\X}$-coresolution of $M$.

(3) If there exists the smallest non-negative integer $n$ such that there is a $\Y_{\X}$-coresolution of length $n$ of $M$, we call ${\rm coresdim}^{\X}_{\Y}(M):=n$ the $\Y_{\X}$-coresolution dimension of $M$. If such $n$ does not exist, we set ${\rm coresdim}^{\X}_{\Y}(M):=\infty$.

(4) We define the $\Y_{\X}$-coresolution dimension of the class $\Z$ as
$${\rm coresdim^{\X}_{\Y}(\Z):=sup\{coresdim^{\X}_{\Y}(Z)|Z\in\Z\}}.$$

(5) We denote by $\Y^{\vee}_{\X,\infty}(\Y^{\vee}_{\X})$ the class of all the objects in $\C$ having a (finite) $\Y_{\X}$-coresolution.

(6) $(\X,\Y)^{\vee}_{\infty}:=\X\cap\Y^{\vee}_{\X,\infty}$ and $(\X,\Y)^{\vee}:=\X\cap\Y^{\vee}_{\X}$.

Dually, we can define the $\Y_{\X}$-resolution, the $\Y_{\X}$-resolution dimension ${\rm resdim}^{\X}_{\Y}(M)$ of $M$, the classes $\Y^{\wedge}_{\X,\infty},\Y^{\wedge}_{\X}$ and $(\Y,\X)^{\wedge}_{\infty}:=\Y^{\wedge}_{\X,\infty}\cap\X , (\Y,\X)^{\wedge}:=\Y^{\wedge}_{\X}\cap\X$.

\end{definition}

By \cite[Lemma 2.13(a)]{OC}, we can easily know the following result.

\begin{lemma}\label{lemma 1}
Let $\C$ be an extriangulated category, $\Y$ and $\X \subseteq\C$. Then ${\rm pd_{\Y}(\X^{\vee})=pd_{\Y}(\X)}$.

\end{lemma}

The following theorem is a generalization of \cite[Theorem 4.4]{MH}, it is the basis of our later conclusions to some extent.

\begin{theorem}\label{th1}
Let $\X\subseteq\Y\subseteq\C$ be both closed under extensions, $\omega$ be a relative generator in $\X$ and $0\in\X$. Then we can get the following statements.

{\rm (1)} For any $Z\in\X^{\vee}_{\Y}$, with ${\rm coresdim}^{\Y}_{\X}(Z):=n$, there are $\mathbb{E}$-triangles
$$\xymatrix{Z\ar[r]^{gz}&M_Z\ar[r]&C_Z\ar@{-->}[r]&} \makebox{with}~~ C_Z\in(\Y,\omega)^{\vee}, M_Z\in\X$$
$$\xymatrix{K_Z\ar[r]&B_Z\ar[r]^{fz}&Z\ar@{-->}[r]&} \makebox{with} ~~B_Z\in\omega^{\vee}_{\Y}, K_Z\in\X,$$
where ${\rm coresdim}^{\Y}_{\omega}(C_Z):=n-1$ and ${\rm coresdim}^{\Y}_{\omega}(B_Z)\leq n$.

{\rm (2)} $B_Z\in(\Y,\omega)^{\vee}$ if $Z\in(\Y,\X)^{\vee}$.

{\rm (3)} Let $\omega\subseteq {^{\perp}\X}$. Then $\omega^{\vee}\subseteq {^{\perp}\X}$, $fz$ is a $\omega^{\vee}$-${\rm precover}$, and $gz$ is an $\X$-${\rm preenvelope}$.
\end{theorem}

\proof (1)We can use inductive hypothesis to prove it since ${\rm coresdim}^{\Y}_{\X}(Z)=n<\infty$.

Let $n=0$. By the definition of $Z\in\X^{\vee}_{\Y}$, we get $Z\in\X\cap\Y$. Because $\omega$ is a relative generator in $\X$, there is an $\mathbb{E}$-triangle $\xymatrix{X'\ar[r]&W\ar[r]&Z\ar@{-->}[r]&}$, with $W\in\omega\subseteq(\Y,\omega)^{\vee},X'\in\X$. This is the second $\mathbb{E}$-triangle that we want to find. And we can easily find the first $\mathbb{E}$-triangle $\xymatrix{Z\ar[r]^{1}&Z\ar[r]&0\ar@{-->}[r]&}$ with $0\in(\Y,\omega)^{\vee}$.

Let $n>0$. Then we have an $\mathbb{E}$-triangle $\xymatrix{Z\ar[r]&X_0\ar[r]^f&Y\ar@{-->}[r]&}$, with $X_0\in\X,Y\in\Y$, where ${\rm coresdim}^{\Y}_{\X}(Y)=n-1$. Then, by inductive hypothesis, we get an $\mathbb{E}$-triangle $\xymatrix{K_Y\ar[r]&B_Y\ar[r]^{f_Y}&Y\ar@{-->}[r]&}$, with $B_Y\in\omega^{\vee}_{\Y}, K_Y\in\X\subseteq\Y$, where ${\rm coresdim}^{\Y}_{\omega}(B_Y)\leq n-1$. By Proposition \ref{prop2},
we have the following commutative diagram
\[
\xy
(-7,21)*+{K_Y}="-12";
(7,21)*+{K_Y}="-14";
(-21,7)*+{Z}="0";
(-7,7)*+{E}="2";
(7,7)*+{B_Y}="4";
(-21,-7)*+{Z}="10";
(-7,-7)*+{X_0}="12";
(7,-7)*+{Y}="14";
{\ar@{=} "-12";"-14"};
{\ar "-12";"2"};
{\ar "-14";"4"};
{\ar "0";"2"};
{\ar "2";"4"};
{\ar@{=} "0";"10"};
{\ar "2";"12"};
{\ar "4";"14"};
{\ar "10";"12"};
{\ar "12";"14"};
{\ar@{}|\circlearrowright "-12";"4"};
{\ar@{}|\circlearrowright "0";"12"};
{\ar@{}|\circlearrowright "2";"14"};
\endxy
\]
of $\E$-triangles. Thus we can get the first $\mathbb{E}$-triangle
$$\xymatrix{Z\ar[r]^{h}&E\ar[r]&B_Y\ar@{-->}[r]&}  \makebox{with} ~E\in\X,B_Y\in\Y.$$
Since $E\in\X$ and $\omega$ is a relative generator in $\X$, there is an $\mathbb{E}$-triangle $\xymatrix{L\ar[r]&W\ar[r]^{h'}&E\ar@{-->}[r]&}$ with $W\in\omega, L\in\X.$

Then, by Proposition \ref{RemET4op}, we have the following commutative diagram
\[
\xy
(-21,7)*+{L}="0";
(-7,7)*+{Z}="2";
(7,7)*+{Z}="4";
(-21,-7)*+{L}="10";
(-7,-7)*+{W}="12";
(7,-7)*+{E}="14";
(-7,-21)*+{B_Y}="22";
(7,-21)*+{B_Y}="24";
{\ar "0";"2"};
{\ar "2";"4"};
{\ar@{=} "0";"10"};
{\ar "2";"12"};
{\ar "4";"14"};
{\ar "10";"12"};
{\ar "12";"14"};
{\ar "12";"22"};
{\ar "14";"24"};
{\ar@{=} "22";"24"};
{\ar@{}|\circlearrowright "0";"12"};
{\ar@{}|\circlearrowright "2";"14"};
{\ar@{}|\circlearrowright "12";"24"};
\endxy
\]
of $\E$-triangles. Thus we can get the second $\mathbb{E}$-triangle
$$\xymatrix{L\ar[r]&Z'\ar[r]&Z\ar@{-->}[r]&}  \makebox{with} ~Z'\in\omega^{\vee}_{\Y},L\in\X.$$
And we find that ${\rm coresdim}^{\Y}_{\omega}(Z')\leq n$.

Since $\omega\subseteq\X$, by the definition of ${\rm coresdim}^{\Y}_{\X}(B_Y)$, we get
$${\rm coresdim^{\Y}_{\X}(B_Y)\leq coresdim^{\Y}_{\omega}(B_Y)\leq n-1}.$$
By the first $\mathbb{E}$-triangle $\xymatrix{Z\ar[r]^{h}&E\ar[r]&B_Y\ar@{-->}[r]&} $, we get
$${\rm \emph{n}=coresdim^{\Y}_{\X}(Z)\leq 1+coresdim^{\Y}_{\X}(B_Y)\leq \emph{n}},$$
It implies that ${\rm coresdim}^{\Y}_{\omega}(B_Y)=n-1$.

(2) We have the $\mathbb{E}$-triangle
$\xymatrix{L\ar[r]&Z'\ar[r]&Z\ar@{-->}[r]&}$, since $\Y$ is closed under extensions and $L,Z\in\Y$, then $Z'\in\Y$.
Hence (2) has been proved.

(3) By Lemma \ref{lemma 1},  we can get $\omega^{\vee}\subseteq ^{\perp}\X$ since $\omega\subseteq ^{\perp}\X$.

Let $Z\in\X^{\vee}_{\Y}$, by the result of (1), we get $C_Z\in ^{\perp}\X$ and $M_Z\in\X\subseteq(\omega^{\vee})^{\perp}$.
For any morphism $f\colon Z\rightarrow X$, with $X\in\X$, we can get the following diagram by the $\mathbb{E}$-triangle:
\[
\xy
(-21,7)*+{Z}="0";
(-7,7)*+{M_Z}="2";
(7,7)*+{C_Z}="4";
(-21,-7)*+{X}="10";
(-7,-7)*+{E}="12";
(7,-7)*+{C_Z}="14";
{\ar^{gz} "0";"2"};
{\ar "2";"4"};
{\ar "0";"10"};
{\ar "2";"12"};
{\ar@{=} "4";"14"};
{\ar "10";"12"};
{\ar "12";"14"};
\endxy
\]
This implies that $gz$ is an $\X$-preenvelope.
\qed

Next we will recall the notion of thick class in an extriangulated category.

\begin{definition}
Let $\C$ be an extriangulated category and $\M,\X\subseteq \C$.

(1) $\M$ is closed under Cones in $\M\cap\X$ if, for any $\mathbb{E}$-triangle
$$\xymatrix{M\ar[r]&M'\ar[r]&M''\ar@{-->}[r]&} \makebox{with}~ M,M'\in \M\cap\X,$$
we have that $M''\in\M$. If $\M\subseteq\X$, we simply say that $\M$ is closed under Cones. And $\M$ is closed under CoCones is defined dually.

(2) $\M$ is right thick if it is closed under extensions, direct summands and Cones.

(3) $\M$ is left thick if it is closed under extensions, direct summands and CoCones.

(4) $\M$ is thick if it is right and left thick.

\end{definition}

The following proposition reveals some connections among $(\X,\Y)^{\vee}_{\infty},(\X,\Y)^{\vee}$ and $(\X,\Y)^{\vee}_{n}$, they have closure properties under certain conditions.

\begin{proposition}\label{prop1}
Let $\C$ be an extriangulated category and $\X,\Y\subseteq \C$, the following statements hold.

{\rm (1)} Let $\Y=\Y^{\oplus<\infty}$ and let $\X$ be closed under extensions and such that $(\X,\Y)^{\vee}_{\infty}\subseteq {^{\perp_1}\Y}$. Then, for a given $\mathbb{E}$-triangle
$\xymatrix{A\ar[r]&B\ar[r]&C\ar@{-->}[r]&}$, with $A,C\in (\X,\Y)^{\vee}_{\infty}$, we have
$${\rm coresdim^{\X}_{\Y}(B)\leq max\{coresdim^{\X}_{\Y}(A),coresdim^{\X}_{\Y}(C)\}}.$$
Furthermore, $(\X,\Y)^{\vee}_{\infty}$ and $(\X,\Y)^{\vee}$ are closed under extensions.

{\rm (2)} Let $\X={\rm smd}(\X)$ and $(\X,\Y)^{\vee}_{\infty}$ be both closed under extensions. Then $(\X,\Y)^{\vee}_{\infty}$ is closed under direct summands.
\end{proposition}

\proof (1) Give an $\mathbb{E}$-triangle $\eta_0: \xymatrix{A\ar[r]^u&B\ar[r]^v&C\ar@{-->}[r]&}$, with $A,C\in (\X,\Y)^{\vee}_{\infty}$. Since $\X$ is closed under extensions and $A,C\in\X$, then we have $B\in\X$. And by definition, we  have $\mathbb{E}$-triangles
$$\eta_A^0: \xymatrix{A\ar[r]^a&Y_A\ar[r]&A_1\ar@{-->}[r]&} \makebox{and}\quad
\eta_C^0: \xymatrix{C\ar[r]^c&Y_C\ar[r]&C_1\ar@{-->}[r]&}$$
with $Y_A,Y_C\in\Y$ and $A_1,C_1\in(\X,\Y)^{\vee}_{\infty}$. Because $C\in(\X,\Y)^{\vee}_{\infty}\subseteq {^{\perp_1}\Y}\subseteq {^{\perp_1}Y_A}$, then by definition, we have the following exact sequence
$$\C(C,Y_A)\rightarrow\C(B,Y_A)\rightarrow\C(A,Y_A)$$
Then we can know that there is a morphism $\alpha:B\rightarrow Y_A$ such that $\alpha u=a$. Consider the morphism $b:=\left(\begin{array}{cc}
\alpha\\
cv
\end{array}\right):B\rightarrow Y_A\oplus Y_C$. Since
$$\left(\begin{array}{cc}
1\\
0
\end{array}\right)a=\left(\begin{array}{cc}
a\\
0
\end{array}\right)=\left(\begin{array}{cc}
\alpha u\\
0
\end{array}\right)=\left(\begin{array}{cc}
\alpha\\
cv
\end{array}\right)u  \quad\makebox{and}\quad \left(\begin{array}{cc}
0&1
\end{array}\right)\left(\begin{array}{cc}
\alpha\\
cv
\end{array}\right)=cv,$$
by \cite[Lemma 3.4]{HZZ}, we can have $\mathbb{E}$-triangles
$$\eta_B^0: \xymatrix{B\ar[r]^b&Y_A\oplus Y_C\ar[r]&B_1\ar@{-->}[r]&} \makebox{and}\quad
\eta_1: \xymatrix{A_1\ar[r]&B_1\ar[r]&C_1\ar@{-->}[r]&}$$
with $A_1,C_1\in(\X,\Y)^{\vee}_{\infty},Y_A\oplus Y_C\in\Y,B_1\in\X$. Then we assume that we have the following $\mathbb{E}$-triangles
$$\eta_B^{k-1}: \xymatrix{B_{k-1}\ar[r]^{b_{k-1}}&Y_{A,k-1}\oplus Y_{C,k-1}\ar[r]&B_{k}\ar@{-->}[r]&} \makebox{and}\quad
\eta_k: \xymatrix{A_k\ar[r]^{u_k}&B_k\ar[r]^{v_k}&C_k\ar@{-->}[r]&}$$
with $A_k,C_k\in(\X,\Y)^{\vee}_{\infty},Y_{A,k-1}\oplus Y_{C,k-1}\in\Y,B_k\in\X,\forall k\leq n$.

By definition, there are $\mathbb{E}$-triangles
$$\eta_A^n: \xymatrix{A_n\ar[r]^{a_n}&Y_{A,n}\ar[r]&A_{n+1}\ar@{-->}[r]&} \makebox{and}\quad
\eta_C^n: \xymatrix{C_n\ar[r]^{c_n}&Y_{C,n}\ar[r]&C_{n+1}\ar@{-->}[r]&}$$
with $Y_{A,n},Y_{C,n}\in\Y$ and $A_{n+1},C_{n+1}\in(\X,\Y)^{\vee}_{\infty}$. Because $C_n\in(\X,\Y)^{\vee}_{\infty}\subseteq {^{\perp_1}\Y}$, then we have the following exact sequence
$$\C(C_n,Y_{A,n})\rightarrow\C(B_n,Y_{A,n})\rightarrow\C(A_n,Y_{A,n})$$
Then we can know that there is a morphism $\alpha_n:B_n\rightarrow Y_{A,n}$ such that $\alpha_n u_n=a_n$. Consider the morphism $b_n:=\left(\begin{array}{cc}
\alpha_n\\
c_nv_n
\end{array}\right):B_n\rightarrow Y_{A,n}\oplus Y_{C,n}$. Since
$$\left(\begin{array}{cc}
1\\
0
\end{array}\right)a_n=\left(\begin{array}{cc}
a_n\\
0
\end{array}\right)=\left(\begin{array}{cc}
\alpha_n u_n\\
0
\end{array}\right)=\left(\begin{array}{cc}
\alpha_n\\
c_nv_n
\end{array}\right)u_n  \quad\makebox{and}\quad \left(\begin{array}{cc}
0&1
\end{array}\right)\left(\begin{array}{cc}
\alpha_n\\
c_nv_n
\end{array}\right)=c_nv_n,$$
by \cite[Lemma 3.4]{HZZ}, we can have $\mathbb{E}$-triangles
$$\eta_B^n: \xymatrix{B_n\ar[r]^{b_n}&Y_{A,n}\oplus Y_{C,n}\ar[r]&B_{n+1}\ar@{-->}[r]&} \makebox{and}\quad
\eta_{n+1}: \xymatrix{A_{n+1}\ar[r]&B_{n+1}\ar[r]&C_{n+1}\ar@{-->}[r]&}$$
with $A_{n+1},C_{n+1}\in(\X,\Y)^{\vee}_{\infty},Y_{A,n}\oplus Y_{C,n}\in\Y,B_{n+1}\in\X$.

By the family of $\mathbb{E}$-triangles $\{\eta_B^i\}_{i=0}^{\infty}$, we can know that $B\in(\X,\Y)^{\vee}_{\infty}$.
Since $A,C\in (\X,\Y)^{\vee}_{\infty}$, then the families of $\mathbb{E}$-triangles $\{\eta_A^k\}_{k=0}^{\infty}$ and $\{\eta_C^k\}_{k=0}^{\infty}$ can form a $(\X,\Y)$-coresolution of minimal length.

For $m:={\rm max}\{{\rm coresdim}^{\X}_{\Y}(A),{\rm coresdim}^{\X}_{\Y}(C)\}$, we get $A_k=0=C_k,\forall k>m$. Consider the family of $\mathbb{E}$-triangles $\{\eta_k\}_{k=1}^{\infty}$, we get $B_k=0,\forall k>m$. So we have ${\rm coresdim}^{\X}_{\Y}(B)\leq m$.

(2) There is an $\mathbb{E}$-triangle $\xymatrix{W\ar[r]&V\ar[r]^f&U\ar@{-->}[r]&}$, with $V\in (\X,\Y)^{\vee}_{\infty}$. Because $\X={\rm smd}(\X)$, we can get $U,W\in\X$. Since $V\in (\X,\Y)^{\vee}_{\infty}$, by definition, there is an $\mathbb{E}$-triangle
$$\xymatrix{V\ar[r]^g&Y_0\ar[r]&V_1\ar@{-->}[r]&}, \makebox{with}~ Y_0\in\Y ~\makebox{and}~ V_1\in (\X,\Y)^{\vee}_{\infty}.$$
By \cite[Definition 2.12(ET4)]{NP}, there exists the following commutative diagram
\[
\xy
(-21,7)*+{W}="0";
(-7,7)*+{V}="2";
(7,7)*+{U}="4";
(-21,-7)*+{W}="10";
(-7,-7)*+{Y_0}="12";
(7,-7)*+{W_1}="14";
(-7,-21)*+{V_1}="22";
(7,-21)*+{V_1}="24";
{\ar "0";"2"};
{\ar "2";"4"};
{\ar@{=} "0";"10"};
{\ar "2";"12"};
{\ar "4";"14"};
{\ar "10";"12"};
{\ar "12";"14"};
{\ar "12";"22"};
{\ar "14";"24"};
{\ar@{=} "22";"24"};
{\ar@{}|\circlearrowright "0";"12"};
{\ar@{}|\circlearrowright "2";"14"};
{\ar@{}|\circlearrowright "12";"24"};
\endxy
\]
of $\E$-triangles. Thus we get the $\mathbb{E}$-triangles $$\eta:\xymatrix{U\ar[r]&W_1\ar[r]&V_1\ar@{-->}[r]&}$$
$$\mu_0:\xymatrix{W\ar[r]&Y_0\ar[r]&W_1\ar@{-->}[r]&}$$
Because $U,V_1\in\X$ and $\X$ is closed under extension, we have $W_1\in\X$. Consider $\mathbb{E}$-triangles $\eta$ and $\xymatrix{W\ar[r]&W\ar[r]&0\ar@{-->}[r]&}$, we can get $\mathbb{E}$-triangle $$\xymatrix{V\ar[r]&W\oplus W_1\ar[r]&V_1\ar@{-->}[r]&}$$
Since $V,V_1\in(\X,\Y)^{\vee}_{\infty}$ and $(\X,\Y)^{\vee}_{\infty}$ is closed under extension, we have $W\oplus W_1\in(\X,\Y)^{\vee}_{\infty}$. Hence we can get a family of $\mathbb{E}$-triangles $\{\mu_i\}^{\infty}_{i=0}$ with $W_0=W,Y_i\in\Y,W_i\in\X,\forall i\geq 0$. So we can get $W\in(\X,\Y)^{\vee}_{\infty}$, it implies $(\X,\Y)^{\vee}_{\infty}$ is closed under direct summands.
\qed
\medskip

We give a necessary and sufficient conditions for some classes to be thick.

\begin{theorem}
Let $(\C,\E,\s)$ be an extriangulated category and $(\X,\Y)\subseteq \C^2$ be satisfied $\Y=\Y^{\oplus <\infty}$, $\X= {\rm smd}(\X)$ is closed under extensions and $\E(\X,\X\cap\Y)=0$. Then the following statements hold.

{\rm (1)} $(\X,\Y)^{\vee}_{\infty}=(\X,\X\cap\Y)^{\vee}_{\infty}$ and it is closed under extensions and direct summands. Further we can know $(\X,\Y)^{\vee}=(\X,\X\cap\Y)^{\vee}$ and it is closed under extensions.

{\rm (2)} $(\X,\Y)^{\vee}_{\infty}$ is left thick if $\X$ is left thick.

\proof (1) Since $\X$ is closed under extensions and direct summands, by definition we can get $(\X,\Y)^{\vee}_{\infty}=(\X,\X\cap\Y)^{\vee}_{\infty}$. Because ${\rm Ext}_{\C}^1(\X,\X\cap\Y)=0$, we have that $(\X,\X\cap\Y)^{\vee}_{\infty}\subseteq\X\subseteq {^{\perp_1}}(\X\cap\Y)$. Hence, by applying Proposition \ref{prop1} (1)(2) to the pair $(\X,\X\cap\Y)$, we can get the result that we want.

(2) Since $\X$ is left thick, then $\X$ is closed under CoCones. Next we will prove that $(\X,\X\cap\Y)^{\vee}_{\infty}$ is closed under CoCones. Consider an $\mathbb{E}$-triangle $$\xymatrix{A\ar[r]^a&B\ar[r]&C\ar@{-->}[r]&}, \makebox{with} ~B,C\in (\X,\X\cap\Y)^{\vee}_{\infty}$$
then there is an $\mathbb{E}$-triangle $$\xymatrix{B\ar[r]^b&W_0\ar[r]&C_0\ar@{-->}[r]&}, \makebox{with} ~W_0\in\X\cap\Y,C_0\in(\X,\X\cap\Y)^{\vee}_{\infty}$$
By \cite[Definition 2.12(ET4)]{NP}, there exists a commutative diagram
\[
\xy
(-21,7)*+{A}="0";
(-7,7)*+{B}="2";
(7,7)*+{C}="4";
(-21,-7)*+{A}="10";
(-7,-7)*+{W_0}="12";
(7,-7)*+{C'}="14";
(-7,-21)*+{C_0}="22";
(7,-21)*+{C_0}="24";
{\ar "0";"2"};
{\ar "2";"4"};
{\ar@{=} "0";"10"};
{\ar "2";"12"};
{\ar "4";"14"};
{\ar "10";"12"};
{\ar "12";"14"};
{\ar "12";"22"};
{\ar "14";"24"};
{\ar@{=} "22";"24"};
{\ar@{}|\circlearrowright "0";"12"};
{\ar@{}|\circlearrowright "2";"14"};
{\ar@{}|\circlearrowright "12";"24"};
\endxy
\]
of $\E$-triangles. Thus we get the $\mathbb{E}$-triangles $$\eta:\xymatrix{A\ar[r]&W_0\ar[r]&C'\ar@{-->}[r]&}$$
$$\eta':\xymatrix{C\ar[r]&C'\ar[r]&C_0\ar@{-->}[r]&}$$
Because $C,C_0\in(\X,\X\cap\Y)^{\vee}_{\infty}$, by (1) $(\X,\X\cap\Y)^{\vee}_{\infty}$ is closed under extensions, then $C'\in(\X,\X\cap\Y)^{\vee}_{\infty}$. Since $W_0\in\X\cap\Y$ and $\X$ is left thick, we get $A\in\X$. Hence $A\in(\X,\X\cap\Y)^{\vee}_{\infty}$ and $(\X,\X\cap\Y)^{\vee}_{\infty}$ is closed under CoCones. Because $(\X,\X\cap\Y)^{\vee}_{\infty}=(\X,\Y)^{\vee}_{\infty}$, $(\X,\Y)^{\vee}_{\infty}$ is left thick.

\qed
\end{theorem}

For the convenience of later discussion, we give the notion of
the relative $n$-quotients ($n$-subobjects) and its relative properties.

\begin{definition}
Let $\C$ be an extriangulated category, $\Y\subseteq\X\subseteq\C$ and $n\geq 1$.

(1) $\Y$ is closed by $n$-quotients in $\X$ if for any complex in $\C$ $$\xymatrix{A\ar[r]&Y_n\ar[r]^{\varphi_n}&Y_{n-1}\ar[r]&\cdots\ar[r]&Y_1\ar[r]^{\varphi_1}&B}$$ with $Y_i\in\Y,{\rm CoCone}(\varphi_i)\in\X,  \forall i\in[1,n]$ and $B\in\X$, we get that $B\in\Y$.

(2) $\Y$ is closed by $n$-subobjects in $\X$ if for any complex in $\C$ $$\xymatrix{A\ar[r]&Y_1\ar[r]^{\varphi_1}&Y_2\ar[r]&\cdots\ar[r]&Y_n\ar[r]^{\varphi_n}&B}$$ with $Y_i\in\Y,{\rm Cone}(\varphi_i)\in\X , \forall i\in[1,n]$ and $A\in\X$, we get that $A\in\Y$.
\end{definition}

\begin{definition}
For an extriangulated category $\C$, $n\geq 1,\X,\T\subseteq \C$. The ${\rm Fac}^{\X}_{n}(\T)$ is defined as the class of all the objects $C\in\C$ satisfying a complex
$$\xymatrix{K\ar[r]&T_n\ar[r]^{f_n}&T_{n-1}\ar[r]&\cdots\ar[r]&T_{2}\ar[r]^{f_2}&T_{1}\ar[r]^{f_1}&C}$$
with ${\rm CoCone}(f_i)\in \X$ and $T_i\in\T\cap\X , \forall i\in [1,n]$.
\end{definition}

\begin{proposition}\label{prop4}
Let $(\C,\E,\s)$ be an extriangulated category, $\X,\T\subseteq\C$ and $\alpha\subseteq\T^{\perp}\cap\X^{\perp}$ be a relative cogenerator in $\X$. Then we have that $\X\cap\T^{\perp}$ is closed by $n$-quotients in $\X$ if and only if ${\rm pd}_{\X}(\T)\leq n,n\geq 1$.
\end{proposition}

\proof ($\Rightarrow$) Let $X\in\X$ and $M\in\T$. Because $\alpha$ is a relative cogenerator in $\X$, by definition, there is a complex
$$\xymatrix{X\ar[r]&I_0\ar[r]&\cdots\ar[r]&Y_{n-1}\ar[r]&V}$$ with $K={\rm CoCone}(f)\in\X,V\in\X$ and $I_i\in\alpha, \forall i\in[0,n-1].$
We can easily get an $\E$-triangle $\xymatrix{K\ar[r]&I_{n-1}\ar[r]&V\ar@{-->}[r]&}$, then we apply the functor ${\rm Hom}_{\C}(M,-)$ to the $\E$-triangle. Because $M\in\T,I_i\in\alpha\subseteq\T^{\perp}\cap\X^{\perp}$, we can get $\E^i(M,I_{n-1})=0=\E^{i+1}(M,I_{n-1})$, it follows that $\E^i(M,V)\cong \E^{i+1}(M,K),\forall i\geq 1$.
By the dual of Shifting Lemma, we get $\E^{i+1}(M,K)\cong \E^{i+n}(M,X)$. Since $\X\cap\T^{\perp}$ is closed by $n$-quotients in $\X$ and $V\in\X\cap\T^{\perp}$, we have
$$0=\E^{i}(M,V)\cong \E^{i+1}(M,K)\cong \E^{i+n}(M,X), \forall i\geq 1.$$
Thus ${\rm pd}_{\X}(M)\leq n$.

($\Leftarrow$) Let $\xymatrix{A\ar[r]&X_n\ar[r]^{\varphi_n}&X_{n-1}\ar[r]&\cdots\ar[r]&X_{1}\ar[r]^{\varphi_1}&B}$ be a complex, with $X_1,\cdots,X_n\in\X\cap\T^{\perp},B\in\X,{\rm CoCone}(\varphi_i)\in\X, \forall i\in[1,n]$. By applying the dual of Shifting Lemma, we get $\E^{k}(M,B)\cong \E^{k+n}(M,A), \forall k\geq 1$, with $\E^{k+n}(M,A)=0 ,\forall k\geq 1$. We know that $A\in\X,{\rm pd}_{\X}(\T)\leq n$, then $\E^{k}(M,B)=0$, therefore $B\in\X\cap\T^{\perp}$.
\qed

\medskip

Next, we will give the notion of $n$-$\X$-tilting class in an extriangulated category. Our goal is to get a generalization of the relative tilting theory on $\X\subseteq \C$.

\begin{definition}\label{definition 1}
Let $\C$ be an extriangulated category and $\X\subseteq \C,n\in\mathbb{N}$. A class $\T\subseteq \C$ is called $n-\X$-tilting if the following conditions hold true.

(T0): $\T={\rm smd}(\T)$.

(T1): ${\rm pd}_{\X}(\T)\leq n$.

(T2): $\T\cap\X\subseteq\T^{\perp}$.

(T3): There is a class $\omega\subseteq\T^{\vee}_{\X}$ which is a relative generator in $\X$.

(T4): There is a class $\alpha\subseteq\X^{\perp}\cap\T^{\perp}$ which is a relative cogenerator in $\X$.

(T5): Every $Z\in\X\cap\T^{\perp}$ admits a $\T$-precover $T'\rightarrow Z$, with $T'\in\X$.

\end{definition}

\begin{lemma}\label{lemma2}
Let $\C$ be an extriangulated category and $\X,\T\subseteq \C$, the following conditions hold.

{\rm (1)} $\T^{\vee}_{\X}\cap\X\subseteq\T^{\vee}\subseteq {^{\perp}(\T^{\perp})}\subseteq {^{\perp}(\T^{\perp}\cap\X)}$.

{\rm (2)} If $\X={\rm smd}(\X)$ and $\T$ satisfies {\rm (T0)} and  {\rm (T2)}, then $(\T\cap\X)^{\vee}_{\X}\cap\T^{\perp}=\T\cap\X$.
\end{lemma}

\proof (1) Applying Lemma \ref{lemma 1}, we know that ${\rm pd}_{\X}(\T^{\vee})={\rm pd}_{\X}(\T)$, it implies that $\T^{\vee}\subseteq {^{\perp}(\T^{\perp})}$, then we can get (1).

(2) Let $A\in(\T\cap\X)^{\vee}_{\X}\cap\T^{\perp}$. Then there is an $\mathbb{E}$-triangle
$\eta:\xymatrix{A\ar[r]&T_0\ar[r]&A'\ar@{-->}[r]&}$
with $T_0\in\T\cap\X$ and $A'\in(\T\cap\X)^{\vee}_{\X}$, where $A'\in{^{\perp}(\T^{\perp})}$.
Because $A\in\T^{\perp}$, we get $\eta$ splits, then $A\in\T\cap\X$, thus $(\T\cap\X)^{\vee}_{\X}\cap\T^{\perp}\subseteq\T\cap\X$.

And by (T2) and $\T\cap\X\subseteq(\T\cap\X)^{\vee}_{\X}$, we know $\T\cap\X\subseteq(\T\cap\X)^{\vee}_{\X}\cap\T^{\perp}$.
\qed

\medskip

The following Lemma gives an important condition that $\T^{\perp}\cap\X\subseteq {\rm Fac}^{\X}_{1}(\T)$, we are interested in studying its relationship with $n$-$\X$-tilting classes. The following results will be expanded around this point.

\begin{lemma}\label{lemma4}
For an extriangulated category $\C$, $\X\subseteq\C$ is closed under extensions. We can get $\T^{\perp}\cap\X\subseteq {\rm Fac}^{\X}_{1}(\T)$ if $\T\subseteq\C$ satisfies {\rm(T3)}.

\proof Let $A\in\T^{\perp}\cap\X$. Since $\T\subseteq\C$ satisfies (T3), there is a class $\omega\subseteq\T^{\vee}_{\X}$ which is a relative generator in $\X$, then there exists an $\mathbb{E}$-triangle
$$\eta_1:\xymatrix{K\ar[r]&W\ar[r]^{a}&A\ar@{-->}[r]&}$$
with $W\in\omega$ and $K\in\X$. Applying Lemma \ref{lemma2}(a) and (T3), there is an $\mathbb{E}$-triangle
$$\eta_2:\xymatrix{W\ar[r]^{b}&B\ar[r]&C\ar@{-->}[r]&}$$
with $B\in\T$ and $C\in\T^{\vee}\cap\X\subseteq{^{\perp}(\T^{\perp})}\cap\X$. Since $\X$ is closed under extensions, we get $B\in\X$.

Applying the dual of Proposition \ref{RemET4op}, we get the following diagram
\[
\xy
(-21,7)*+{K}="0";
(-7,7)*+{K}="2";
(-21,-7)*+{W}="10";
(-7,-7)*+{B}="12";
(7,-7)*+{C}="14";
(-21,-21)*+{A}="20";
(-7,-21)*+{B'}="22";
(7,-21)*+{C}="24";
{\ar@{=} "0";"2"};
{\ar "0";"10"};
{\ar "2";"12"};
{\ar_{b} "10";"12"};
{\ar "12";"14"};
{\ar_{a} "10";"20"};
{\ar^{x} "12";"22"};
{\ar@{=} "14";"24"};
{\ar^{t} "20";"22"};
{\ar "22";"24"};
{\ar@{}|\circlearrowright "0";"12"};
{\ar@{}|\circlearrowright "10";"22"};
{\ar@{}|\circlearrowright "12";"24"};
\endxy
\]
of $\E$-triangles. Then there exist two $\mathbb{E}$-triangles
$$\eta_3:\xymatrix{K\ar[r]&B\ar[r]^x&B'\ar@{-->}[r]&},\eta_4:\xymatrix{A\ar[r]^t&B'\ar[r]&C\ar@{-->}[r]&}$$
by definition,we have $B'\in {\rm Fac}^{\X}_{1}(\T)$. Since $A\in\T^{\perp},C\in{^{\perp}(\T^{\perp})}$, then $\eta_4$ splits.
There exists $y:B'\rightarrow A$ such that $yt=1_A$. By \cite[Proposition 3.17]{NP}, and the $\mathbb{E}$-triangle
$$\eta_5:\xymatrix{K'\ar[r]&B\ar[r]^{yx}&A\ar@{-->}[r]&},$$ we get the following diagram
\[
\xy
(-21,7)*+{K}="0";
(-7,7)*+{W}="2";
(7,7)*+{A}="4";
(-21,-7)*+{K'}="10";
(-7,-7)*+{B}="12";
(7,-7)*+{A}="14";
(-21,-21)*+{C}="20";
(-7,-21)*+{C}="22";
{\ar "0";"2"};
{\ar^{a} "2";"4"};
{\ar "0";"10"};
{\ar^{b} "2";"12"};
{\ar@{=} "4";"14"};
{\ar "10";"12"};
{\ar_{yx} "12";"14"};
{\ar "10";"20"};
{\ar "12";"22"};
{\ar@{=} "20";"22"};
{\ar@{}|\circlearrowright "0";"12"};
{\ar@{}|\circlearrowright "2";"14"};
{\ar@{}|\circlearrowright "10";"22"};
\endxy
\]
Then we can get an $\mathbb{E}$-triangle
$\xymatrix{K\ar[r]&K'\ar[r]&C\ar@{-->}[r]&}$.
Since $\X$ is closed under extensions and $K,C\in\X$, then $K'\in\X$. Therefore $\eta_5$ satisfies the definition of ${\rm Fac}^{\X}_{1}(\T)$, we get $A\in {\rm Fac}^{\X}_{1}(\T)$.
\qed

\end{lemma}

\begin{lemma}\label{lemma3}
Let $\C$ be an extriangulated category, $\X={\rm smd}(\X)\subseteq\C$ is closed under extensions, $\T\subseteq\C$ satisfies {\rm (T2), (T5)} and such that $\T^{\perp}\cap\X\subseteq {\rm Fac}^{\X}_{1}(\T)$, the following conditions hold.

{\rm (1)} $\T\cap\X$ is a relative generator in $\T^{\perp}\cap\X$.

{\rm (2)} For each morphism $A\rightarrow X$, where $A\in{^{\perp}(\T^{\perp}\cap\X)}$ and $X\in\T^{\perp}\cap\X$, factors through $\T\cap\X$. And if $\T={\rm smd}(\T)$, then we have
$$\T\cap\X=\T^{\perp}\cap\X\cap{^{\perp}(\T^{\perp}\cap\X)}=\T^{\perp}\cap\X\cap{^{\perp}(\T^{\perp})}.$$
\end{lemma}

\proof (1) We know $\T\cap\X\subseteq\T^{\perp}$ from (T2). Let $X\in\T^{\perp}\cap\X$, Then there admits a $\T$-precover $g:T'\rightarrow X$ with $T'\in\X$ by (T5). Assume there exists an $\mathbb{E}$-triangle
$\xymatrix{K\ar[r]&T'\ar[r]^g&X\ar@{-->}[r]&}$, we can get $K\in\T^{\perp}$ since g is a $\T$-precover and $\T\subseteq\C$ satisfies (T2). By the definition of $X\in {\rm Fac}^{\X}_{1}(\T)$, there exists an $\mathbb{E}$-triangle
$$\xymatrix{K'\ar[r]&B\ar[r]^f&X\ar@{-->}[r]&}$$
with $B\in \T\cap\X$ and $K'\in\X$. Applying Proposition \ref{prop2}, we have the following diagram
\[
\xy
(-7,21)*+{K'}="-12";
(7,21)*+{K'}="-14";
(-21,7)*+{K}="0";
(-7,7)*+{Z}="2";
(7,7)*+{B}="4";
(-21,-7)*+{K}="10";
(-7,-7)*+{T'}="12";
(7,-7)*+{X}="14";
{\ar@{=} "-12";"-14"};
{\ar "-12";"2"};
{\ar "-14";"4"};
{\ar "0";"2"};
{\ar "2";"4"};
{\ar@{=} "0";"10"};
{\ar "2";"12"};
{\ar^f "4";"14"};
{\ar "10";"12"};
{\ar^g "12";"14"};
{\ar@{}|\circlearrowright "-12";"4"};
{\ar@{}|\circlearrowright "0";"12"};
{\ar@{}|\circlearrowright "2";"14"};
\endxy
\]
Then there are two $\mathbb{E}$-triangles
$$\eta_1:\xymatrix{K'\ar[r]&Z\ar[r]&T'\ar@{-->}[r]&},\eta_2:\xymatrix{K\ar[r]&Z\ar[r]&B\ar@{-->}[r]&}$$
Since $\X={\rm smd}(\X)\subseteq\C$ is closed under extensions and $K',T'\in\X$, then $Z\in\X$. And $\eta_2$ splits since $K\in T^{\perp},B\in\T$, we get $K\in\X$. Hence $K\in \T^{\perp}\cap\X$, $\T\cap\X$ is a relative generator in $\T^{\perp}\cap\X$.

(2) Assume there is a morphism $f:A\rightarrow X$ with $A\in{^{\perp}(\T^{\perp}\cap\X)}$ and $X\in\T^{\perp}\cap\X$. We can know there exists an $\mathbb{E}$-triangle
$$\xymatrix{L\ar[r]^{k}&T'\ar[r]^g&X\ar@{-->}[r]&}$$
with $T'\in\T\cap\X$ and $L\in\T^{\perp}\cap\X$. Obviously $\E(A,L)=0$, then we get ${\rm Hom}_{\C}(A,g)$ is surjective and $f$ factors through $g$. Note that for every $A'\in\T^{\perp}\cap\X\cap{^{\perp}(\T^{\perp})}~ (\rm or  ~ \T^{\perp}\cap\X\cap{^{\perp}(\T^{\perp}\cap\X)})$, we get $1_{A'}$ factors through $\T\cap\X$.
\qed

\begin{lemma}\label{lemma5}
For an extriangulated category $\C$,$\X={\rm smd}(\X)\subseteq\C$ is closed under extensions, let $\T={\rm smd}(\T)\subseteq\C$ satisfy {\rm (T1), (T2), (T4), (T5)} and such that $\T^{\perp}\cap\X\subseteq {\rm Fac}^{\X}_1(\T)$. We also get $\X\subseteq(\T^{\perp}\cap\X)^{\vee}_{\X}$ and $(\T\cap\X)^{\vee}\subseteq {^{\perp}(\T^{\perp}\cap\X)}$, and for each $X\in\X$, the following conditions hold.

{\rm (1)} $ m:={\rm coresdim}^{\X}_{\T^{\perp}\cap\X}(X)\leq {\rm pd}_{\X}(\T)<\infty;$

{\rm (2)} there exist two $\mathbb{E}$-triangles
$$\xymatrix{X\ar[r]&M_X\ar[r]&C_X\ar@{-->}[r]&}  \makebox{and}~\xymatrix{K_X\ar[r]&B_X\ar[r]&X\ar@{-->}[r]&}$$
with $M_X,K_X\in\T^{\perp}\cap\X;C_X,B_X\in\X; {\rm coresdim}^{\X}_{\T\cap\X}(C_X)=m-1$ and ${\rm coresdim}^{\X}_{\T\cap\X}(B_X)\leq m$;

{\rm (3)} $B_X\rightarrow X$ is a $(\T\cap\X)^{\vee}$-${\rm precover}$;

{\rm (4)} $X\rightarrow M_X$ is a $\T^{\perp}\cap\X$-${\rm preenvelope}$.
\end{lemma}

\proof We can easily get $\T\cap\X$ is a $\T^{\perp}\cap\X$-projective relative generator in $\T^{\perp}\cap\X$ and $\X\subseteq(\T^{\perp}\cap\X)^{\vee}_{\X}$ by Lemma \ref{lemma3}. Then, we can prove the lemma by applying Theorem \ref{th1}.
\qed

\begin{proposition}\label{prop3}
For an extriangulated category $\C$, $\X={\rm smd}(\X)\subseteq\C$ is closed under extensions, let $\T={\rm smd}(\T)\subseteq\C$ satisfy {\rm (T1),(T2),(T4),(T5)}. We get that T satisfies {\rm (T3)} if and only if $\T^{\perp}\cap\X\subseteq {\rm Fac}^{\X}_1(\T)$ and we can choose a relative generator $\omega$ in $\X$ such that $\omega\subseteq(\T\cap\X)^{\vee}_{\X}$.
\end{proposition}

\proof
$(\Rightarrow)$  We can prove that if $\T\subseteq\C$ satisfies (T3). It follows that $\T^{\perp}\cap\X\subseteq {\rm Fac}^{\X}_{1}(\T)$ by Lemma \ref{lemma4}.

$(\Leftarrow)$ By Lemma \ref{lemma5}, for each $X\in\X$, there is an $\mathbb{E}$-triangle
$$\xymatrix{X\ar[r]&M_X\ar[r]&C_X\ar@{-->}[r]&},$$
with $M_X\in\T^{\perp}\cap\X$ and $C_X\in(\X,\T\cap\X)^{\vee}_{\X}$.
By the definition of $\T^{\perp}\cap\X\subseteq {\rm Fac}^{\X}_1(\T)$, we can get $M_X$ admits an $\mathbb{E}$-triangle
$$\xymatrix{M_X'\ar[r]&T_0\ar[r]&M_X\ar@{-->}[r]&},$$
with $T_0\in\T\cap\X$ and $M_X'\in\X$. By the axiom (ET4$^{\rm op}$), we have the following diagram
\[
\xy
(-21,7)*+{M_X'}="0";
(-7,7)*+{M_X'}="2";
(-21,-7)*+{P_X}="10";
(-7,-7)*+{T_0}="12";
(7,-7)*+{C_X}="14";
(-21,-21)*+{X}="20";
(-7,-21)*+{M_X}="22";
(7,-21)*+{C_X}="24";
{\ar@{=} "0";"2"};
{\ar "0";"10"};
{\ar "2";"12"};
{\ar "10";"12"};
{\ar "12";"14"};
{\ar "10";"20"};
{\ar "12";"22"};
{\ar@{=} "14";"24"};
{\ar "20";"22"};
{\ar "22";"24"};
{\ar@{}|\circlearrowright "0";"12"};
{\ar@{}|\circlearrowright "10";"22"};
{\ar@{}|\circlearrowright "12";"24"};
\endxy
\]
Then there exists an $\mathbb{E}$-triangle
$$\xymatrix{M_X'\ar[r]&P_X\ar[r]&X\ar@{-->}[r]&}~ \makebox{with} ~P_X\in(\X,\T\cap\X)^{\vee}_{\X}.$$
So $P_X$ is a relative generator in $\X$ satisfying (T3).
\qed

\medskip

We recall the notions of $\X$-hereditary and $\X$-complete pair for the discussion of our main result.

\begin{definition}
For an extriangulated category $\C$ and $\X\subseteq\C$, a pair $(\A,\B)\subseteq\C^2$ is called $\X$-hereditary if ${\rm id}_{\A\cap\X}(\B\cap\X)=0$. If $\X=\C$, we call $(\A,\B)$ is hereditary.
\end{definition}

\begin{definition}
Let $\C$ be an extriangulated category, $(\A,\B)\subseteq\C^2$. We say $(\A,\B)$ is left $\X$-complete if for each $X\in\X$, there is an $\mathbb{E}$-triangle $\xymatrix{B\ar[r]&A\ar[r]&X\ar@{-->}[r]&}$, with $A\in\A\cap\X$ and $B\in\B\cap\X$. The right $\X$-complete pair is defined dually. Moreover, $(\A,\B)$ is called $\X$-complete if it is left and right $\X$-complete.
\end{definition}

The following result will enrich the tilting theory.

\begin{theorem}\label{th2}
Let $(\C,\E,\s)$ be an extriangulated category, $\X={\rm smd}(\X)\subseteq\C$ be closed under extensions and $\T\subseteq\C$ satisfying {\rm (T1), (T2), (T3), (T4) and (T5)}, the following conditions hold.

{\rm (1)} ${^{\perp}(\T^{\perp}\cap\X)}\cap \X={^{\perp}(\T^{\perp})}\cap \X=\T^{\vee}_{\X}\cap\X=(\T\cap\X)^{\vee}_{\X}\cap\X $ if $\T={\rm smd}(\T)$.

{\rm (2)} $\T^{\perp}\cap\X={\rm Fac}^{\X}_k(\T)\cap\X, \forall k\geq {\rm max}\{1,{\rm pd}_{\X}(\T)\}.$

{\rm (3)} If $\T={\rm smd}(\T)$, then $({^{\perp}(\T^{\perp})},\T^{\perp})$ is $\X$-complete and hereditary.

\proof
(1) Applying Lemma \ref{lemma2}(1), $\T^{\vee}_{\X}\cap\X\subseteq\T^{\vee}\subseteq {^{\perp}(\T^{\perp})}\subseteq {^{\perp}(\T^{\perp}\cap\X)}$, we can get that
$$(\T\cap\X)^{\vee}_{\X}\cap\X \subseteq\T^{\vee}_{\X}\cap\X\subseteq{^{\perp}(\T^{\perp})}\cap\X\subseteq {^{\perp}(\T^{\perp}\cap\X)}\cap\X.$$
To prove (1), we firstly consider $X\in{^{\perp}(\T^{\perp}\cap\X)}\cap\X$. We can know that there is an $\mathbb{E}$-triangle by Lemma \ref{lemma5} and Proposition \ref{prop3}
 $$\xymatrix{X\ar[r]&M_X\ar[r]&C_X\ar@{-->}[r]&}$$
with $M_X\in\T^{\perp}\cap\X$ and $C_X\in(\X,\T\cap\X)^{\vee}_{\X}$.
By Lemma \ref{lemma2}(1), we also get $C_X\in{^{\perp}(\T^{\perp}\cap\X)}$. Since ${^{\perp}(\T^{\perp}\cap\X)}$ is closed under extensions and $X,C_X\in{^{\perp}(\T^{\perp}\cap\X)}$, then $M_X\in{^{\perp}(\T^{\perp}\cap\X)}$.
We know that $M_X\in\T\cap\X$ since $\T\cap\X=\T^{\perp}\cap\X\cap{^{\perp}(\T^{\perp}\cap\X)}$. Hence, by definition, $X\in(\T\cap\X)^{\vee}_{\X}\cap\X$.

(2) We know that $\T\cap\X$ is a relative generator in $\T^{\perp}\cap\X$ by Lemma \ref{lemma4} and Lemma \ref{lemma3}(1). We also get $\T^{\perp}\cap\X\subseteq {\rm Fac}^{\X}_k(\T)\cap\X ,\forall k\geq 1$.

Assume $m:={\rm max}\{1,{\rm pd}_{\X}(\T)\}$, let $C\in {\rm Fac}^{\X}_k(\T)\cap\X$ with $k\geq m$. Then there is a complex
$$\xymatrix{K\ar[r]&T_k\ar[r]^{f_k}&\cdots\ar[r]^f_2&T_{1}\ar[r]^{f_1}&C}$$
with ${\rm CoCone}(f_i)\in \X$ and $T_i\in\T\cap\X ,\forall i\in [1,k]$. By (T1),(T2),(T4)and Proposition \ref{prop4} $\X\cap\T^{\perp}$ is closed by n-quotients in $\X$, we get $C\in\X\cap\T^{\perp}$.

(3) By definition, we can easily get the pair $({^{\perp}(\T^{\perp})},\T^{\perp})$ is hereditary. We know that
$(\T\cap\X)^{\vee}_{\X}\cap\X \subseteq{^{\perp}(\T^{\perp})}\cap\X$ by Lemma \ref{lemma2}(1).

Applying Lemma \ref{lemma5} and Proposition \ref{prop3}, for every $X\in\X$, there exist two $\mathbb{E}$-triangles
$$\xymatrix{X\ar[r]&M_X\ar[r]&C_X\ar@{-->}[r]&}  \makebox{and}~\xymatrix{K_X\ar[r]&B_X\ar[r]&X\ar@{-->}[r]&}$$
with $M_X,K_X\in\T^{\perp}\cap\X;C_X,B_X\in (\T\cap\X)^{\vee}_{\X}\cap\X$.
It satisfies the definition of $\X$-complete.
\qed

\end{theorem}

On the basis of Theorem \ref{th2}, we can get the following result.

\begin{theorem}
Let $(\C,\E,\s)$ be an extriangulated category, $n\geq 1,\X={\rm smd}(\X)\subseteq\C$ be closed under extensions and $\T={\rm smd}(\T)\subseteq\C$ satisfying {\rm (T4)} and {\rm (T5)}, the following conditions are equivalent.

{\rm (1)} $\T$ is $n$-$\X$-${\rm tilting}$.

{\rm (2)} ${\rm Fac}^{\X}_n(\T)\cap\X=\T^{\perp}\cap\X$.

{\rm (3)} $\T^{\perp}\cap\X={\rm Fac}^{\X}_k(\T)\cap\X ,\forall k\geq n$.

{\rm (4)} $\T^{\perp}\cap\X$ is closed by $n$-quotients in $\X$ and $\T\cap\X\subseteq\T^{\perp}\cap\X\subseteq {\rm Fac}^{\X}_1(\T)$.

\proof

$(1)\Rightarrow (2)$  It can be proved by Theorem \ref{th2}(2).

$(2)\Rightarrow (3)$  It is obvious that ${\rm Fac}^{\X}_{n+1}(\T)\cap\X\subseteq {\rm Fac}^{\X}_n(\T)\cap\X$. Next we will prove that ${\rm Fac}^{\X}_{n+1}(\T)\cap\X\supseteq {\rm Fac}^{\X}_n(\T)\cap\X$.

Assume $N\in {\rm Fac}^{\X}_n(\T)\cap\X=\T^{\perp}\cap\X$, since $\T$ satisfies (T5), there is a $\T$-precover $f:A\rightarrow N$ with $A\in\X$. Because ${\rm Fac}^{\X}_{n}(\T)\subseteq {\rm Fac}^{\X}_1(\T)$, we get that $f$ is a deflation. Then there exists an $\mathbb{E}$-triangle
$$\eta:\xymatrix{K\ar[r]&A\ar[r]^f&N\ar@{-->}[r]&}$$
with $A\in\T\cap\X\subseteq {\rm Fac}^{\X}_n(\T)\cap\X=\T^{\perp}\cap\X$. We get that $K\in\T^{\perp}$ by $A,N\in\T^{\perp}$ and $f$ is a $\T$-precover.

Next we prove that $K\in\X$. Since $N\in {\rm Fac}^{\X}_1(\T)$, by definition, there is an $\mathbb{E}$-triangle
$$\eta_1:\xymatrix{K'\ar[r]&M_0\ar[r]^{f'}&N\ar@{-->}[r]&}$$
with $M_0\in\T\cap\X$ and $K'\in\X$. By Proposition \ref{prop2}, we have the following diagram
\[
\xy
(-7,21)*+{K'}="-12";
(7,21)*+{K'}="-14";
(-21,7)*+{K}="0";
(-7,7)*+{P}="2";
(7,7)*+{M_0}="4";
(-21,-7)*+{K}="10";
(-7,-7)*+{A}="12";
(7,-7)*+{N}="14";
{\ar@{=} "-12";"-14"};
{\ar "-12";"2"};
{\ar "-14";"4"};
{\ar "0";"2"};
{\ar "2";"4"};
{\ar@{=} "0";"10"};
{\ar "2";"12"};
{\ar "4";"14"};
{\ar "10";"12"};
{\ar "12";"14"};
{\ar@{}|\circlearrowright "-12";"4"};
{\ar@{}|\circlearrowright "0";"12"};
{\ar@{}|\circlearrowright "2";"14"};
\endxy
\]
Then we get an an $\mathbb{E}$-triangle
$$\eta_2:\xymatrix{K\ar[r]&P\ar[r]&M_0\ar@{-->}[r]&} \makebox{and}~~ \eta_3:\xymatrix{K'\ar[r]&P\ar[r]&A\ar@{-->}[r]&}$$
Because $\X={\rm smd}(\X)\subseteq\C$ is closed under extensions and $K',A\in \X$, then $P\in\X$. Since $M_0\in\T\cap\X$ and $K\in\T^{\perp}$, then $\eta_2$ splits and $K\in\X$. Hence, we have $K\in\T^{\perp}\cap\X={\rm Fac}^{\X}_n(\T)\cap\X$, by definition, $N\in {\rm Fac}^{\X}_{n+1}(\T)\cap\X$.

$(3)\Rightarrow (4)$
By (3),  we have $\T\cap\X\subseteq {\rm Fac}^{\X}_{n}(\T)\cap\X=\T^{\perp}\cap\X\subseteq {\rm Fac}^{\X}_{1}(\T)$. Since $$\T^{\perp}\cap\X={\rm Fac}^{\X}_{n}(\T)\cap\X={\rm Fac}^{\X}_{n+1}(\T)\cap\X,$$ we can get that $\T^{\perp}\cap\X$ is closed by $n$-quotients in $\X$.

$(4)\Rightarrow (1)$
By (4), (T4) and Proposition \ref{prop4}, we have that $${\rm pd}_{\X}(\T)\leq n,\T^{\perp}\cap\X={\rm Fac}^{\X}_1(\T)\cap\X$$ and (T1), (T2) hold. Hence, by
applying Proposition \ref{prop3}, we can prove (1).
\qed

\end{theorem}

\medskip

\hspace{-5mm}\textbf{Chenbei Xie}\\[1mm]
College of Mathematics, Hunan Institute of Science and Technology, 414006 Yueyang, Hunan, P. R. China\\[1mm]
E-mail: xcb19980313@163.com


\begin{thebibliography}{9}

\bibitem{AB} M. Auslander, R. Buchweitz.
The homological theory of maximal Cohen-Macaulay approximations. (French summary)
Colloque en l'honneur de Pierre Samuel (Orsay, 1987).
M\'{e}m. Soc. Math. France (N.S.) No. 38 (1989), 5-37.


\bibitem{GNP}  M. Gorsky,  H. Nakaoka,  Y. Palu. Positive and negative extensions in extriangulated categories. arXiv: 2103.12482, 2021.

\bibitem{H}
M. Hashimoto.
Auslander-Buchweitz approximations of equivariant modules.
London Mathematical Society Lecture Note Series, 282. Cambridge University Press, Cambridge, 2000.

\bibitem{HZZ}  J. Hu,  D. Zhang,  P. Zhou. Proper resolutions and Gorensteinness in extriangulated categories. Frontiers of Mathematics in China. Front. Math. China 16 (2021), no. 1, 95-117.

\bibitem{MH}  A. A. Monroy, O. M. Hern\'{a}ndez. Relative tilting theory in abelian categories I: Auslander-Buchweitz-Reiten approximations theory in subcategories and cotorsion pairs.
    arXiv: 2104.11361, 2021.

\bibitem{MH2}  A. A. Monroy, O. M. Hern\'{a}ndez. Relative tilting theory in abelian categories II:
$n$-$\X$-tilting theory. arXiv: 2112.14873, 2021.

\bibitem{MSSS} O. Mendoza Hern\'{a}ndez, E. S\'{a}enz Valadez, V. Santiago Vargas, M. Souto Salorio.
Auslander-Buchweitz approximation theory for triangulated categories.
Appl. Categ. Structures 21 (2013), no. 2, 119-139.



\bibitem{NP} H. Nakaoka, Y. Palu.  . Cah. Topol. G\'{e}om. Diff\'{e}r. Cat\'{e}g. 60 (2019), no. 2, 117-193.


\bibitem{OC} C. Mendoza, C. S\'{a}enz.
Tilting categories with applications to stratifying systems. J. Algebra 302 (2006), no. 1, 419-449.
\end{thebibliography}
\end{document}